\documentclass[a4paper]{amsart}\usepackage{amsfonts,amssymb,bbm,mathrsfs}
\usepackage{enumerate}

\def\C{\Bbb{C}}\def\k{\mathbbm{k}}\def\K{\mathbb{K}}

\def\N{\Bbb{N}}

\newcommand{\quots}[2]{{\footnotesize\left.\raisebox{0.4ex}{$#1$}\! / \!\raisebox{-0.4ex}{$#2$}\right.}}

\def\tM{\tilde{M}}\def\tm{\tilde{m}}
\def\tn{\tilde{n}}\def\tq{\tilde{q}}

\def\tw{{\tilde{w}}}

\def\hM{{\hat{M}}}\def\hR{\widehat{R}}

\def\al{\alpha}\def\be{\beta}\def\de{\delta}
\def\ep{\epsilon}

\def\ca{\mathfrak a}
\def\cC{\mathcal C}

\def\cO{\mathcal O}\def\cP{\mathscr P}
\def\cm{{\frak m}}\def\cp{{\frak p}}

\def\uf{{\underline{f}}}\def\ug{{\underline{g}}}\def\uh{\underline{h}}
\def\uq{\underline{q}}
\def\ux{\underline{x}}\def\uy{{\underline{y}}}

\newcommand{\ber}{\begin{array}{l}}\newcommand{\eer}{\end{array}}
\newcommand{\bpm}{\begin{pmatrix}}\newcommand{\epm}{\end{pmatrix}}
\newcommand{\bM}{\begin{matrix}}\newcommand{\eM}{\end{matrix}}
\newcommand{\bee}{\begin{enumerate}}\newcommand{\eee}{\end{enumerate}}
\newcommand{\bei}{\begin{itemize}}\newcommand{\eei}{\end{itemize}}

\def\sset{\subset}\def\sseteq{\subseteq}\def\smin{\setminus}

\def\Mat{Mat_{m\times n}(R)}

\newtheorem{Lemma}{Lemma}[section]\newcommand{\bel}{\begin{Lemma}}\newcommand{\eel}{\end{Lemma}}

\newtheorem{Theorem}[Lemma]{Theorem}\newcommand{\bthe}{\begin{Theorem}}\newcommand{\ethe}{\end{Theorem}}
\newtheorem{Proposition}[Lemma]{Proposition}\newcommand{\bprop}{\begin{Proposition}}\newcommand{\eprop}{\end{Proposition}}
\newtheorem{Corollary}[Lemma]{Corollary}\newcommand{\bcor}{\begin{Corollary}}\newcommand{\ecor}{\end{Corollary}}
\newtheorem{Definition}[Lemma]{Definition}\newcommand{\bed}{\begin{Definition}}\newcommand{\eed}{\end{Definition}}
\newtheorem{Definition-Proposition}[Lemma]{Definition-Proposition}

\def\bpr{~\\{\em Proof.\ }}
\newcommand{\epr}{{\hfill\ensuremath\blacksquare}}
\newtheorem{Remark}[Lemma]{Remark}\newcommand{\beR}{\begin{Remark}\rm}\newcommand{\eeR}{\end{Remark}}
\newtheorem{Example}[Lemma]{Example}\newcommand{\bex}{\begin{Example}\rm}\newcommand{\eex}{\end{Example}}
\newtheorem{Problem}[Lemma]{Problem}\newcommand{\bprob}{\begin{Problem}\rm}\newcommand{\eprob}{\end{Problem}}

\newcommand{\bet}{\begin{tabular}{cccccccc}}\newcommand{\eet}{\end{tabular}}
\newcommand{\beq}{\begin{equation}}\newcommand{\eeq}{\end{equation}}

\newcommand{\bin}[2]{\binom{#1}{#2}}

\title[]{S\MakeLowercase{ome genericity results over} N\MakeLowercase{oetherian rings}}
\author{D\MakeLowercase{mitry} K\MakeLowercase{erner}}
\address{Department of Mathematics, Ben Gurion University of the Negev, P.O.B. 653, Be'er Sheva 84105, Israel.}
\email{dmitry.kerner@gmail.com}
\date{\today
}
\thanks{D.K. was supported by Israel Science Foundation, grant No. 844/14.}
\thanks{We thank L. Avramov,  G.M. Greuel and A.F. Boix for important advices}

\begin{document}\maketitle
\begin{abstract}
Let $M$ be a filtered module.
Some properties of elements of $M$ are ``generic" in the following sense:
\bei
\item (being open/stable) if an element $z\in M$ has a property $\cP$ then any approximation of $z$ has $\cP$;
\item  (being dense)  any element of $M$ is approximated by an element that has $\cP$.
\eei
(Here the approximation is taken in the filtered sense.)

Moreover, one can often ensure an approximation  with further special properties, e.g. avoiding a prescribed set of submodules.

We prove that being a regular sequence is a generic property. As immediate applications we get corollaries on
 the generic grades of modules, heights of ideals, properties of determinantal ideals, acyclicity of generalized Eagon-Northcott complexes and vanishing of Tor/Ext.
\end{abstract}
\setcounter{secnumdepth}{6} \setcounter{tocdepth}{2}

\section{Introduction}

\subsection{}
Given an $R$-module $M$  with a filtration $M=M_0\supsetneq M_1\supsetneq\cdots$, and a certain property $\cP$, denote by $M^\cP\sseteq M$ the subset of elements
 that possess $\cP$.
\bed\label{Def.Genericity.Property}
A property $\cP$ is said to hold $M_\bullet$-generically if $M^\cP$ is open and dense in the $M_\bullet$-topology.
\eed
Explicitly this means:
\bei
\item (being open) for any $z\in M^\cP$ and some $n\gg1$ (that depends on $z$) holds: $\{z\}+M_n\sseteq M^\cP$.
\item (being dense) for any $z\in M$ and any $n\in \N$ there exists $w\in M_n$ such that $z+w\in M^\cP$.

Equivalently, for any $n\in \N$, the projection $M\stackrel{\pi_n}{\to}\quots{M}{M_n}$ satisfies:   $\pi_n(M)=\pi_n(M^\cP)$.
\eei

A property is often ``expected"  if it holds over the ring of indeterminates, $R[\ux]$ or $R[[\ux]]$.
 But such an expected property can be violated by a (degenerate) specialization, $R[\ux]\to R$, $R[[\ux]]\to R$, $\ux\to \uf$.
 We prove that many expected properties hold generically on the module of specializations.

\subsection{Conventions}\label{Sec.Conventions}
The rings in this paper are (unital, commutative) Noetherian, with no assumptions on characteristic.
 All the modules are finitely-generated.
 Unless stated otherwise $\cm\sset R$ denotes the   Jacobson radical. (In particular, for a local ring, $\cm$ is the maximal ideal.)

By $L,M,N$ we denote $R$-modules.

We use the multi-indices, $\ux=(x_1,\dots,x_r)$, $\uf=(f_1,\dots,f_r)$, etc.

\subsection{The content of the paper}
\bee[{\bf \S 1}] \setcounter{enumi}{1}
\item 
addresses the semi-continuity/stability of grade, exactness and height.
 In lemma \ref{Thm.Stability.of.exactness.under.small.deformations} we recall that exactness of a complex is preserved under ``higher order perturbations",
  for the filtration $\{\cm^n\}$.
 Thus exactness is an open property for $\cm$-adic filtrations.
 This implies, e.g. semi-continuity of the grade of ideal and stability of a sequence being regular.
 In lemma \ref{Thm.Stability.Height.under.small.deformations} we prove that the height of ideals is semi-continuous.

For non-local rings the Jacobson radicals are often small (sometimes zero). Thus one would like to replace the filtration $\{\cm^n\}$ by $\{I^n\}$,
 for some other ideal $I\sset R$.
Our lemmas are sharp, if the elements of a filtration are not contained in $\{\cm^n\}$ then exactness
 is not preserved, see remark \ref{Rem.why.Jacboson.radical}.

Most of these results are well known
  (see the references in \S\ref{Sec.Stability.Semicontinuity}),  we give short proofs for completeness.

\

\item 
is about approximations with expected grade/height/exactness.
Given an $M$-regular sequence of polynomials/power series, $\uh(\uf):=(h_1(\ux),\dots,h_c(\ux))\in R[\ux]^c,R[[\ux]]^c$, its specialization, $\uh(\uf)\in R^c$,
  is not necessarily $M$-regular.
In theorem \ref{Thm.Deforming.to.expected.grade} we prove that this specialization is approximated by a regular specialization,
 $\uh(\uf+\ug)\in R^c$.
 In particular, the ideal $(\uh(\uf))\sset R$ is approximated  by a complete intersection ideal, of grade $c$, as expected.
 Here the approximation $\ug$ can be chosen small in a very strict sense and can avoid a prescribed set of submodules.

Note that this approximation, $(\uh(\uf))\rightsquigarrow (\uh(\uf+\ug))$, is not a flat deformation.

 Similarly, in theorem \ref{Thm.Deforming.to.expected.height} we prove: if $height(\uh)=c$ (even if this sequence is non-regular) then
  the ideal $(\uh(\uf))$ is approximated by an ideal $(\uh(\uf+\ug))$ of height $c$.

Thus the set of ``good" specializations is dense.

Among examples and applications of this are:
 \bei
\item   sequences are approximated by  regular sequences, ideals are approximated by complete intersections;
\item any (scheme) subgerm   is approximated by a subgerm of expected codimension.
\eei

\item
consists of the applications of these stability/approximation results to ideals/modules/complexes. We obtain characterization of
 ``How does the generic object look like?".
\bei
\item In \ref{Sec.Application.to.Matrices} we establish the ``generic" properties of determinantal ideals/schemes:
any matrix is approximated by a  matrix whose determinantal ideals are of expected grades/heights. If $(R,\cm)$ is a local Cohen-Macaulay
 ring then any matrix is approximated by one whose determinantal ideals are perfect Cohen-Macaulay.

\item In \S\ref{Sec.Depth.Grade.Improving} we show that the behaviour of depth/grade under projections/specializations
 (which is pathological in the degenerate case) can be improved to the ``expected" one.
 \item In \S\ref{Sec.Application.to.E.N.complexes} we obtain:
 any morphism of modules is approximated by a morphism whose generalized Eagon-Northcott complexes are acyclic.
\item In \S\ref{Sec.Application.to.Tor.Ext} we obtain the generic vanishing of $Tor$ and $Ext$.
\eei

\eee

Some of these statements are well known in Commutative Algebra, some others are known in Singularity Theory for the particular case, $R=\C\{\ux\}$.

 But we could  not find the references for this description of the ``generic object",  in the generality of (local) Noetherian rings.

\subsection{Continuation}
Recall the classical Thom's transversality theorem:  any $C^\infty$-map of     manifolds, $X\stackrel{}{\to}Y$,
 can be made transversal to a given submanifold, $Y\supset Z$,  by arbitrarily small deformation.

Our initial motivation was to convert this into a purely algebraic statement, with purely algebraic proofs.
 (Instead of the classical  proofs using the real/complex topology.)

In the subsequent paper we address the genericity properties of irreducibility, primeness, normality, (non-)regular/singular loci.

\section{Stability/opennes/semi-continuity properties of exactness, grade and height}\label{Sec.Stability.Semicontinuity}
\subsection{Deformations of high order preserve exactness of a complex} Let $(R,\cm)$ be as in \S\ref{Sec.Conventions}
\bel\label{Thm.Stability.of.exactness.under.small.deformations}
Given an exact complex of (f.g.) modules, $(M_\bullet,\phi_\bullet)$.
\bee[1.]
\item  
There exists a sequence of integers $\{n_i\}$ satisfying:
  if  $\{\psi_i\in Hom(M_i,\cm^{n_{i-1}} M_{i-1})\}_i$ and
 the chain of maps, $(M_\bullet,\phi_\bullet+\psi_\bullet)$, is again a complex,  then
 this complex is exact.
\item
 For any prime ideal $\cp\sset R$ there exists a sequence of integers $\{n_i\}$ satisfying:
  if  $\{\psi_i\in Hom(M_i,\cp^{n_{i-1}} M_{i-1})\}_i$ and
 the $\cp$-localized chain of maps, $(M_\bullet,\phi_\bullet+\psi_\bullet)_\cp$, is again a complex,  then
 this complex is exact.
 \eee
\eel

If a complex is bounded then we can choose one uniform integer, $n$.

This statement is well known, see e.g. Proposition 1.1 of   \cite{Eisenbud-Huneke} and page 62 of \cite{Eisenbud.74}.

\bpr {\bf Part 1.}
  It is enough to check the exactness of the deformed complex, $(M_\bullet,\phi_\bullet+\psi_\bullet)$, at each place separately.
 To simplify the notations we start with  an exact sequence $A\stackrel{\phi}{\to}B\stackrel{\psi}{\to}C$
  and establish exactness of the deformation $A\stackrel{\phi+\psi_\phi}{\to}B\stackrel{\psi+\psi_\psi}{\to}C$.

Let $b\in ker(\psi+\psi_\psi)$ thus $\psi(b)=-\psi_\psi(b)\in \cm^n\cdot C$. Thus $\psi(b)\in Im(\psi)\cap (\cm^n\cdot C)\sset Im(\psi|_{\cm^{n-\tn}\cdot C})$,
  by Artin-Rees.
 Hence $\psi(b)=\psi(\de_b)$, for some $\de_b\in \cm^{n-\tn}\cdot B$, and thus $b-\de_b\in Im(\phi)$. Put $b-\de_b=\phi(a)$, then
   $b-(\phi(a)+\psi_\phi(a))\in Ker(\psi+\psi_\psi)\cap (\cm^{n-\tn}\cdot B)$.
    Therefore
\beq
Ker(\psi+\psi_\psi)\sseteq Im(\phi+\psi_\phi)+Ker(\psi+\psi_\psi)\cap (\cm^{n-\tn}\cdot B).
\eeq
Iterate this procedure, to get (by Artin-Rees lemma): $Ker(\psi+\psi_\psi)\sseteq Im(\phi+\psi_\phi)+\cm^n\cdot Ker(\psi+\psi_\psi)$
 for any $n$. Finally, apply the Nakayama lemma for $(R,\cm)$.

{\bf Part 2.} follows now by localization at $\cp$.
\epr

\bex\label{Ex.Preserving.Exactness.Grade.under.deformations} Let $\cm\sset R$ be  as in \S\ref{Sec.Conventions}.
\bee[\bf i.]
\item
If $\uf=(f_1,\dots,f_r)$ is a regular sequence then exists $n\in \N$ such that for any $\ug=(g_1,\dots,g_r)\in \cm^n\cdot R^r$  the sequence
 $\uf+\ug$ is regular. (Apply the lemma to the Koszul complex $K_\uf$:  $\cdots\wedge^2 R^r\stackrel{\uf}{\to} R^r\stackrel{\uf}{\to} R$.)
In particular, if $f\in R$  is not a zero divisor then $f+g$ is not a zero divisor, for any $g\in \cm^n$ with $n\gg1$.

\item
Geometrically, given a (algebraic/analytic/formal) germ $(X,o)$, with its local ring $\cO_{(X,o)}$,
  if $V(\uf)\sset (X,o)$ is a complete intersection then $V(\uf+\ug)\sset (X,o)$ is a complete intersection
 for any $\ug\in \cm^n\cdot \cO_{(X,o)}^r$, with $n\gg1$.
\item More generally, given a module $M$, if a tuple $\uf\in R^r$ is $M$-regular then $\uf+\ug$ is $M$-regular too, for any $\ug\in \cm^n\cdot R^r$ with $n\gg1$.
 Indeed, one applies the lemma to the Koszul complex $K_\uf\otimes M$.
 (This gives Corollary 1.2 of \cite{Huneke-Trivedi}.) 

Equivalently, if the Koszul cohomology vanishes, $H^{i>0}(\uf;M)=0$, then its approximation vanishes too,  $H^{i>0}(\uf+\ug;M)=0$.
\eee

\beR\label{Rem.why.Jacboson.radical}
Recall that for non-local rings the Jacobson radical can be rather small, e.g. $J(\k[\ux])=0$. Therefore in
 part 1 of lemma \ref{Thm.Stability.of.exactness.under.small.deformations} one would like to replace $\cm$ by some ideal $I$ and to weaken
  the assumptions at least to:

{\em For any  exact complex of (f.g.) modules, $(M_\bullet,\phi_\bullet)$, with $\{\phi_i\in Hom(M_i,I^{k_{i-1}} M_{i-1})\}_i$,
 there exists a sequence of integers $\{n_i\gg k_i\}$ satisfying:
  if  $\{\psi_i\in Hom(M_i,I^{n_{i-1}} M_{i-1})\}_i$ and
 the chain of maps  $(M_\bullet,\phi_\bullet+\psi_\bullet)$ is again a complex, then
 this complex is exact.}

   This is obstructed by the   following examples:
\bee[i.]
\item
Let $R$ be a non-local ring, with  a maximal ideal $\ca\sset R$. Assume that $R$ contains a copy of the field $\k:=\quots{R}{\ca}$.
 For $x\in R$ take the exact complex $R\stackrel{x}{\to}R\to\quots{R}{(x)}$.
 Suppose an ideal $I\sset R$ satisfies: $I^n\not\sseteq \ca$ for any $n$.
 Take the filtration $\{I^n\}$.
 Fix some $f\in I^n\smin (I^n\cap \ca)$ and denote the image of $f$ in $\quots{R}{\ca}$ by $c\in \k$. Using the embedding $\k\sset R$ identify $c\in R$.
  Thus $\frac{f}{c}\in I^n\smin (I^n\cap \ca)$ goes to $1\in \k$. But then the complex  $R\stackrel{x-\frac{f}{c}x}{\to}R\to\quots{R}{(x)}$ cannot be exact.
 And one can choose $x\in I^j$, $j\gg1$, and   any $n$.

\item Let $R=\quots{\k[x_1,x_2]}{(x_1\cdot x_2)}$ and $I=(1+x_1)$. Let $f=1+x_1$ (a regular element) and $g=(1+x_1)^n$. For any $n\gg1$ the element $f-g$ is not regular.

 \eee
\eeR

\subsection{Semi-continuity of grade/depth}\label{Sec.Depth.Grade.Semi-continuity}
Fix an ideal $I\sset R$, an $R$-module $M$ and a  sequence $\uf\in \cm\cdot R^r$.
\bee[\bf i.]
\item There exists $n\in \N$ such that  for any  sequence $\ug\in \cm^n\cdot R^r$,
 holds: $grade_{(\uf)}M\le grade_{(\uf+\ug)}M$.
\bpr Choose a maximal $M$-regular subsequence among $\uf$. As $\uf$ is a sequence in the Jacobson radical $\cm\sset R$,
 the length of this maximal regular subsequence equals $grade_{(\uf)}M$.
  Apply part iii. of example \ref{Ex.Preserving.Exactness.Grade.under.deformations}.\epr

Remark that this semi-continuity does not hold if one varies the module. For example, let $R=\quots{\k[[x,\ep]]}{(\ep^2)}$ and $M=coker(R\stackrel{\ep}{\to}R)$.
 Then $grade_\cm M=1$, but the grade of $\tM:=coker(R\stackrel{\ep+x^n}{\to}R)$ is zero, for any $n\in \N$.

\item Take  the  projection $R\stackrel{\pi_\uf}{\to}\quots{R}{(\uf)}$.
 We get the ideal $\pi_\uf I\sset \quots{R}{(\uf)}$ and the module $\pi_\uf M:=M\otimes_R \quots{R}{(\uf)}$.
 Suppose $\uf\in \cm\cdot R^r$ is $M$-regular sequence  then (for the given $I,M$) exists $n\in \N$ such that for any $\ug\in \cm^n\cdot R^r$ holds:
$grade_{\pi_\uf I }\pi_\uf M\le grade_{\pi_{\uf+\ug} I}\pi_{\uf+\ug} M$.
\bpr
Let $\{\tq_\al\}$ be a maximal $\pi_\uf M$-regular sequence in $\pi_\uf I$. Fix some preimages $\{q_\al\}$ in $I$, then $\uf,\uq$ is an $M$-regular sequence.
 Thus, for $n\gg1$, $\uf+\ug,\uq$ is $M$-regular as well. Hence $\pi_{\uf+\ug}(\uq)$ is $\pi_{\uf+\ug}M$-regular.
\epr

The $M$-regularity assumption on $\uf$ is  important. For example, let $r=1$ and $I=\cm$. Take $0\neq f\in Ann(M)$ and $g$ such that $f+g$ is $M$-regular. Then
 $grade_{\pi_f I }\pi_f M= grade_{\pi_{f+g} I}\pi_{f+g} M+1$.

On the other hand, we do not assume the sequence $\uf$ to be $R$-regular.


We list two particular cases:
\bei
\item For $M=R$ we get: the grade of the projected ideal does not decrease under small deformations of the projections.
\item A particular case of projection is the specialization. Suppose $S[\ux]\sseteq R\sseteq S[[\ux]]$, fix some $\uf\in \cm\cdot S^r$ and take the
 specialization $R\to \quots{R}{(\ux-\uf)}$. Accordingly we have:  $I|_\uf:=\quots{I}{(\ux-\uf)}$
 and $M|_\uf:=\quots{M}{(\ux-\uf)M}$.
Then there  exists $n\gg1$ such that for any $\ug\in \cm^n\cdot R^r$, holds:
\[
grade_{I|_{\uf+\ug}}M|_{\uf+\ug}\ge grade_{I|_{\uf}}M|_{\uf}.
\]
\eei

\item
Fix a finite set of triples, $\{\uf^{(\al)},I_\al,M_\al\}_\al$, as in part ii. By repeating part ii., there exists $n\in \N$ such that for any sequences
 $\{\ug^{(\al)}\in \cm^n\cdot R^r\}_\al$ holds:
 \beq
\forall\ \al:\quad grade_{\pi_{\uf^{(\al)}}I_\al}\pi_{\uf^{(\al)}}M_\al\le grade_{\pi_{\uf^{(\al)}+\ug^{(\al)}}I_\al}\pi_{\uf^{(\al)}+\ug^{(\al)}}M_\al.
\eeq
\eee
\eex

\subsection{Semi-continuity of height of ideals}
\bel\label{Thm.Stability.Height.under.small.deformations}
Take  a tuple of elements  $\uf\in R^r$. There exists $n\in \N$ such that for any $\ug\in \cm^n\cdot R^r$ holds:
\[height(\uf+\ug)\ge height(\uf).
\]
\eel
\bpr
Suppose $height(\uf)=l$. Choose a minimal set of elements $x_1,\dots,x_l\in \cm$ such that the ideal  $\sum (x_i)+\sum (f_j)$ is $\cm$-primary.
 Fix $n$ such that  $\sum (x_i)+\sum (f_j)\supseteq\cm^{n-1}$. Then for any $\ug\in \cm^n\cdot R^r$ holds (by Nakayama):
 \beq
\sum (x_i)+\sum (f_j+g_j)=\sum (x_i)+\sum (f_j)\sset R. \quad\quad\quad\quad\quad
 \eeq
Thus $height(\uf+\ug)\ge height(\uf)$. \epr

This statement does not hold for non-local rings without further restrictions, see remark \ref{Rem.why.Jacboson.radical}.ii.

\subsection{The geometric interpretation}
 For a local ring $(R,\cm)$ take the corresponding germ of scheme, $(X,o)=Spec(R)$.
 We have the subgerms,
 $V(\uf)$, $\{V(I_\al)\}_\al\sset (X,o)$. Take their intersections, $\{V(\uf)\cap V(I_\al)\}_\al$.
 For any $\ug\in \cm^n\cdot R^r$, with $n\gg1$, holds:
\beq
\forall\ \al:\quad \quad codim\big(V(\uf+\ug)\cap V(I_\al)\big)\ge  codim\big(V(\uf)\cap V(I_\al)\big).
\eeq
Here the intersections can be not pure co-dimensional, then we take the smallest among the co-dimensions of the irreducible components.

\section{Approximations  with expected grade/height/exactness}\label{Sec.Approximation.expected.grade.exactness}
Given a sequence of polynomials in indeterminates, $\uh\in R[\ux]^c$, one often specializes/restricts to $\uh(\uf)\in  R^c$,
 for some tuple  $\uf\in  R^r$. More generally, for a sequence  of power series in indeterminates, $\uh\in R[[\ux]]^c$,
   and a tuple   $\uf\in \cm\cdot R^r$, one specializes  to $\uh(\uf)\in (\hR^{(\cm)})^c$.
    Here the completion $\hR^{(\cm)}$ is taken to ensure the convergence.

 Even if the sequence $\uh$ is regular, the specialization is not necessarily regular.
 Thus the grade/height of ideals/modules can drop down under specialization:
 \beq
 grade_{(\uh(\uf))}(\hR^{(\cm)})|le grade_{(\uh)}(R[[\ux]]),\quad\quad\quad\quad
 height(\uh(\uf))\le height(\uh).
\eeq
 In this section we prove that  the specializations that preserve regularity/exactness/grade/height
  are dense in    the sense of definition \ref{Def.Genericity.Property}.

\subsection{Improving the regularity of a sequence}
Let $(R,\cm)$, $M$ be as in \S\ref{Sec.Conventions}. Define $M[x]:=M\otimes_R R[\ux]$, resp. $M[[\ux]]:=M\otimes_R R[[\ux]]$.

\bthe\label{Thm.Deforming.to.expected.grade}
Suppose $grade_\cm(M)\ge c$. Suppose a tuple of polynomials, $\uh\in R[\ux]^c$, is $M[\ux]$-regular.
  (Suppose a tuple of power series, $\uh\in R[[\ux]]^c$, is $M[[\ux]]$-regular.)
 Take  a tuple   $\uf\in   R^r$ (for polynomials), resp.   $\uf\in \cm \cdot R^r$ (for power series).
 Take a submodule $N\sseteq \cm\cdot R^r$ such that the annihilator ideal $Ann(\quots{R^r}{N})$ contains an $M$-regular sequence of length $c$.
\bee[1.]
\item  There exists $\ug\in N$ such that the sequence $\uh(\uf+\ug)=\{h_1(\uf+\ug),\dots,h_c(\uf+\ug)\}$
 is   $M$-regular   (for polynomials), resp. an $\hM^{(\cm)}$-regular sequence   (for power series).

\item  Moreover, one can ensure that both $\ug$ and $\uf+\ug$ avoid any given finite collection
   of submodules $\{N_\al\sset R^r\}_\al$ that satisfy: $N_\al\not\supseteq \cm^n\cdot N$  for any  $n,\al$.
\\If $R$ is $\cm$-complete then   one can avoid a countable collection of such submodules  $\{N_\al\sset R^r\}$.
\eee
\ethe
\bpr Part 1. is proved in Steps 1-2. Part 2. is proved in Step 3.

\bee[\bf Step 1.]
\item(Reduction to the particular case: $c=1$ and $M=R$, a domain) Recall that
  $\uq=(q_1,\dots,q_c)$ is $M$-regular iff $q_{i+1}$ is regular on $\quots{M}{(q_1,\dots,q_i)M}$,
for any $i=0,\dots,c-1$.

In addition, as the sequence $\uh(\ux)$ is $M[\ux]$-regular, for any $i=1,\dots,c$ holds:
\[
grade\quots{M[\ux]}{(h_1(\ux),\dots,h_i(\ux))M[\ux]}=grade(M[\ux])-i\ge c+r-i.
\]
(And similarly for $M[[\ux]]$.)

Therefore, by induction, it is enough to verify just the case of one element, $c=1$.

 Note that $N$ does not change during the induction, thus for each $i$
  the ideal  $Ann(\quots{R^r}{N})$ contains an $\quots{M}{(q_1,\dots,q_i)M}$-regular sequence of length $c-i$.

\

Therefore we start with an $M$-regular element $h(\ux)$, i.e. $h(\ux)\not\in Ass(M[x])$, and a submodule $N\sset \cm\cdot R^r$, such that
 the annihilator $Ann\quots{R^r}{N}$ contains an $M$-regular element. For any tuple $\uf\in R^r$ we should find a tuple $\ug\in N$ such that
 $h(\uf+\ug)\not\in Ass(M)$. (For power series:  $h(\uf+\ug)\not\in Ass(\hM^{(\cm)})$.)

By Noetherianity $Ass(M)$ is a finite set of prime ideals in $R$. (Note also: $Ass(M[\ux])=Ass(M)\otimes_R R[\ux]$, $Ass(M[[\ux]])=Ass(M)\otimes_R R[[\ux]]$.)
 By example \ref{Ex.Preserving.Exactness.Grade.under.deformations}.i
 it is enough to avoid each of these primes separately. So, we should prove: if $h(\ux)\not\in \cp[\ux]$ then   exists $\ug$ such that $h(\uf+\ug)\not\in \cp$.
 Replace $R$ by the domain $\quots{R}{\cp}$ and note the  $M$-regular element of  $Ann\quots{R^r}{N}$ remains $M$-regular.

Therefore it is enough to prove:
 \beq\label{Eq.inside.proof.particular.case.to.prove}
\ber
 \text{If $0\neq h\in R[\ux]$ then, for any $\uf\in R^r$ and any $N\sset R^r$ such $Ann\quots{R^r}{N}$ contains}
 \\\text{ an $M$-regular element, exists $\ug\in N\sset R^r$ satisfying: $h(\uf+\ug)\neq 0$.}
 \eer
 \eeq

(And similarly in the case of power series, with $\hR^{(\cm)}$ a domain.)

\item (The case: $c=1$ and  $R$, resp $\hR^{(\cm)}$, is a domain) We prove the statement \eqref{Eq.inside.proof.particular.case.to.prove}.
\bee[i.]
\item(Polynomial case) Suppose $0\neq h\in R[\ux]$. Pass to the field of fractions,
   $\K:=Frac(R)$, it is infinite as $dim(R)\ge 1$. Thus the hypersurface $V(h)\sset \K^r$ is a proper closed subset, even though $\K$ is not algebraically closed.
    Suppose $h(\uf)=0$,
    i.e. $\uf\in V(h)$. Take any point $\tilde{\ug}\in \K^r\smin V(h)$ and consider the line $\{\uf+t \tilde{\ug}|\ t\in \K\}\sset \K^r$.
     The tuple $\tilde{\ug}$ does not necessarily belong to $R^r$, due to its denominators. Rescale $t$ to clear the denominators and moreover to ensure
 $ \tilde{\ug}\in N\sset R^r\sset \K^r$. (This is possible as $Ann(\quots{R^r}{N})$ contains a non-zero element.)

  Finally, this line does not lie fully inside $V(h)$, thus it intersect this hypersurface in a finite number of points.
 Therefore there exists $t_0$ satisfying:
 $ t_0\cdot\tilde{\ug}\in N$  and   $\uf+t_0\cdot \tilde{\ug}\not\in V(h)$.
  Take $g:=t_0\cdot \tilde{\ug}$.
\item(Power series case) Suppose $0\neq h\in R[[\ux]]$.
 We can assume $\uf=0$, by expanding $h(\ux)$ at $\uf$.
 (No convergence problems as $\hR^{(\cm)}$ is  complete   and $\uf\in \cm $.)

Let $k:=ord_{(\ux)}h(\ux)$, accordingly present $h(\ux)=h_k(\ux)+h_{>k}(\ux)$, where $h_k$ is homogeneous of degree $k$, while $h_{>k}(\ux)\in (\ux)^{k+1}$.
 By part i. there exists $\ug\in N$ such that $h_k(\ug)\neq0$. If $ ord_{\cm } h_k(\ug)<ord_{(\ux)}h_{>k}(\ux)$, then obviously $h(\ug)\neq0$.
 Otherwise choose any general $0\neq t\in \cm $ and observe:
 \beq\ber
 ord_{\cm } h_k(t\cdot\ug)= k\cdot ord_{\cm }(t)+ord_{\cm } h_k(\ug),
 \\
ord_{\cm } h_{>k}(t\cdot\ug)\ge (k+1)\cdot ord_{\cm }(t)+ord_{\cm } h_{>k}(\ug).
 \eer\eeq
Thus, for $n\gg1$ we get: $ ord_{\cm } h_k(t^n\ug)<ord_{(\ux)}h_{>k}(t^n\ug)$. And therefore $h(t^n\ug)\neq0$.
\eee

\item (Proof of part 2 of the theorem)
Suppose the deformed sequence $\uh(\uf+\ug)\in R^c$ is $M$-regular.
 (Resp. $\uh(\uf+\ug)\in(\hR^{(\cm)})^c$ is $\hM^{(\cm)}$-regular.)
Fix a finite or countable set of submodules $\{N_\al\sset R^r\}$, as in the assumption.

An auxiliary observation: for any $z\in N$, and any $n,\al$, there exists $w_{z,n,\al}\in \cm^n\cdot N$ such that $z+w_{n,\al}\not\in N_\al$.
  Indeed, we can assume $z\in N_\al$. Take some $\tw\in \quots{N+N_\al}{N_\al}$ that is not annihilated by any power of $\cm$.
   Choose some preimage $w\in N$ and consider the line $\{z+t\cdot w\}_{t\in R}$.
  Take   $t_0\in \cm^n$ that does not annihilate $\tw$, then $z+t_0\cdot w\not\in N_\al$.

Using this auxiliary observation fix some $\uq^{(1)}\in \cm^{n_1}\cdot N$   such that $\ug+  \uq^{(1)} \not\in N_1$ and
 $\uf+\ug+  \uq^{(1)} \not\in N_1$.
  Deform $\uf+\ug\rightsquigarrow \uf+\ug+ \uq^{(1)}$. For $n_1\gg1$ the deformed sequence $\uh(\uf+\ug+\uq^{(1)})$ remains regular,
   by example \ref{Ex.Preserving.Exactness.Grade.under.deformations}.i.i. Now deform to $\uf+\ug+\uq^{(1)}+\uq^{(2)}$, for $n_2\gg n_1$, and so on.
    Eventually, for a finite collection $\{N_\al\}$,
      we get a regular sequence $\uh(\uf+\ug+\sum  \uq^{(i)})$, that satisfies:
\beq
      \ug+\sum  \uq^{(i)}\not\in N_\al\quad \text{ and }
        \uf+\ug+\sum  \uq^{(i)}\not\in N_\al,\quad \text{ for any } \al.
\eeq

For a countable collection, and $R$ being $\cm$-complete, one takes the infinite sum $\sum^\infty  \uq^{(i)} $.
\epr
\eee

\bex\label{Ex.Deforming.to.reg.sequence} As the simplest submodule we take  $N=\oplus^r_{j=1} I_j\sset R^r$,
 where $I_j\sset \cm$ and $\{grade_{I_j}M\ge c\}_j$.
 Then $Ann\quots{R^r}{N}=\cap I_j$ is of grade $min\{grade_{I_j}M\}$, which is at least $c$.
\bee[i.]
\item In the simplest case, $h_1=x_1$, \dots, $h_r=x_r$, we get: for any sequence $\uf\in \cm \cdot R^r$ and any ideals $\{I_j\sset \cm \}$,
  with $\{grade_{I_j}M\ge r\}$,
 there exist  $\{g_j\in I_j\}$ such that
 the sequence $\uf+\ug$ is $M$-regular.
 For $r=1$ and $M=R$ this gives Theorem 2.1 of \cite{Avramov-Iyengar}.

Thus, when $grade_\cm M\ge r$, being $M$-regular sequence is a generic property for the elements of the module $\cm \cdot R^r$,   in the sense of definition \ref{Def.Genericity.Property}.

\item
 Geometrically, for the case $R$-local,  any subscheme $V(\uf)\sset Spec(R)$ is approximated by a complete intersection.

 \item
An equivalent formulation: if the tensored Koszul complex $K_{\uh(\ux)}\otimes M[\ux]$ is exact, then its specialization,   $K_{\uh(\uf)}\otimes M$,
 is not necessarily exact, but is approximated by an exact one,
  $K_{\uh(\uf+\ug)}\otimes M$.
    \item Suppose the sequence $f_1,\dots,f_r$ contains an $M$-regular  subsequence, $f_1,\dots,f_k$, for some $k\le r$. Then for any ideals
   $I_{k+1},\dots,I_r$ satisfying $grade(I_j+\sum^k_{i=1}(f_i))\ge r$, there exist  $g_{k+1}\in I_{k+1},\dots,g_r\in I_r $ such that the sequence
     \[f_1,\dots,f_k,f_{k+1}+g_{k+1},\dots, f_r+g_r
     \] is $M$-regular.
 (In the theorem choose $h_i=f_i$ for $i=1\dots k$ and $h_i=x_i$ for $i=k+1\dots r$.)

In particular, in theorem \ref{Thm.Deforming.to.expected.grade} we do not assume $r\ge c$, the case $r<c$ is also important.

\eee
\eex
\bex\label{Ex.Deforming.to.reg.sequence.several.sequences}
\bee[i.]
\item In the (finite) direct sum case, $M=\oplus M_\al$, the theorem ensures $\ug\in N$ such that the sequence $\uh(\uf+\ug)$
 is simultaneously regular for all $\{M_\al\}$.
 This follows also from example \ref{Ex.Preserving.Exactness.Grade.under.deformations}.ii. Moreover, that example implies: if $R$ is complete,
  then one can ensure the simultaneous $\{M_\al\}$-regularity of $\uh(\uf+\ug)$ for a countable collection of modules.
\item
Theorem \ref{Thm.Deforming.to.expected.grade} can be applied repetitively  to deform several sequences to regular sequences, as follows.
 Take some regular sequences, $\{\uh^{(j)}\in R[[\ux]]^{c}\}_{j=1\dots p}$  and ideals $\{I_j\}_{j=1\dots p}$ with $grade_{I_j}(R)\ge c$. Fix some $r$-tuple
  $\uf\in \cm \cdot R^r$. Fix $\ug^{(1)}\in I_1\cdot R^r$ such that the sequence $\uh^{(1)}(\uf+\ug^{(1)})$ is regular.
   Now
   shrink the ideals $I_2,\dots,I_p$ to the ideals  $I_2\cap \cm ^n,\dots,I_p\cap \cm ^n$, $n\gg1$, so that deformation by any
     $\ug^{(j)}\in I_j\cap \cm ^n$ preserves the regularity of $\uh^{(1)}(\uf+\ug^{(1)})$.
 (This is possible by   example \ref{Ex.Preserving.Exactness.Grade.under.deformations}.)
     Fix $\ug^{(2)}\in (I_2\cap \cm ^n)\cdot R^r$ to ensure the regularity
     of both
 $\uh^{(1)}(\uf+\ug^{(1)}+\ug^{(2)})$ and  $\uh^{(2)}(\uf+\ug^{(1)}+\ug^{(2)})$. And so on.
 \eee
\eex

\beR\label{Rem.Cannot.weaken.assumptions} (Sharpness/Weakening the assumptions)
\bee[i.]
\item The assumption ``$Ann\quots{R^r}{N}$ contains an $M$-regular sequence of length $c$" is needed because of the trivial examples: $c=r$,
 $\uh(\ux)=\ux$, and all the permutations.
\item One would like to weaken the assumptions on $\{I_j\}$ in example \ref{Ex.Deforming.to.reg.sequence}  to:
\[
\{grade_{I_j}M\ge1\}_j\quad \quad\text{ and }\quad\quad grade_{\sum I_j}M\ge c.
\]
 These assumptions do not suffice due to the following example. Choose $\uh\in R[[\ux]]$ as $h_1,\dots,h_{c-1}\in I_1$, while $h_c=x_1$. Suppose $f=0$.
 Then to deform $\uh(\uf)$ to a regular sequence one needs: $grade_{I_1}M\ge c$.

\item
 Theorem \ref{Thm.Deforming.to.expected.grade} addresses the rings $R[\ux]$, $R[[\ux]]$. Geometrically one has
  the subscheme of the affine space $V(\uh)\sset \k^r$ or $V(\uh)\sset (\k^r,o)$
   (a formal germs), intersected with the image of $\uf$. And one deforms inside  $\k^r$ or
   $(\k^r,o)$. One would like to deform also inside other smooth varieties/germs. Namely, one would like to replace
  $R[\ux]$, $R[[\ux]]$ by a regular $R$-algebra  $S$. The statement of the theorem
   holds indeed for some particular algebras, e.g. $R=\k[\uy]_{(\uy)}$, $S=\k[\ux,\uy]_{(\ux,\uy)}$.
     But in general (even if $S$ is local) the statement
  is obstructed by the following geometric reasons. Suppose $Spec(R)$ is the germ of an affine smooth rational curve, while $Spec(S)=Spec(R)\times(C_{g>0},o)$,
  here $(C_{g>0},o)$ is the germ of a smooth non-rational curve. Then the only possible map $Spec(R)\to Spec(S)$ has  the constant second coordinate.
   To deform this map to a proper intersection with the constant section one should pass
    to \'{e}tale  neighborhoods, i.e. henselizations of $R$, $S$.

    \item  In part 2 of theorem \ref{Thm.Deforming.to.expected.grade} one would like to weaken the assumption ``$R$ is $\cm$-complete".
For example,  one might ask: ``$R$ is $\cm$-local and henselian". This weakening is obstructed by the following example.
Let $c=1=r$ and
$M=R=\k\langle\uy\rangle$, $N=\cm\cdot R$, where
  the field $\k$ is at most countable, and $dim(R)\ge2$. Take the set of all the principal primes, $\{\cp_\al\}$, it is countable.
   (As the map sending each algebraic power series to its minimal polynomial sends the set $\{\cp_\al\}$ to the set of all polynomials, $\k[\uy]$.)
   For each $\al$ we have $\sqrt{Ann(\quots{R}{\cp_\al})}\not\supseteq \cm$. And no choice of $g\in \cm$ in part 2. can avoid this set.

\eee
\eeR
\bex\label{Ex.Deforming.complexes}
\bee[i.]
\item
Fix a morphism of free modules, $R^a\stackrel{\phi}{\to}R^b$, $a\le b$. For any ideal $I\subsetneq R$ that contains a non-zero divisor,
 there exists $\psi\in Hom(R^a,I\cdot R^b)$,
 such that the morphism $R^a\stackrel{\phi+\psi}{\to}R^b$ is injective.

Proof: It is enough to check the case $a=b$, i.e. to check the injectivity for a (maximal) square submatrix of $\phi$. Thus (for $a=b$) it is enough
  to verify: $det(\phi+\psi)\in R$ is not a zero divisor.
   This follows straight from theorem \ref{Thm.Deforming.to.expected.grade}, for $c=1$ and $h(\ux)=$the determinant of matrix of indeterminates.

The assumption  that $I$ contains a non-zero divisor is important here. For example, let $R=\quots{\k[[\ux,\uy]]}{(\ux)\cdot(\ux,\uy)}$, where $\k$ is a field.
 This ring contains no regular elements. And the morphism $R\stackrel{x+\psi}{\to}R$ is non-injective for any $\psi\in\cm=(\ux,\uy)$.
\item
 One cannot extend part i. to ``a  complex is approximated by an exact one", as a non-surjective morphism is not approximated by a surjective one.
 Indeed,   $M\stackrel{\phi}{\to}N$ is a non-surjective morphism of modules iff its  fiber at the origin,
   $\quots{M}{\cm\cdot M}\stackrel{\phi\otimes\quots{R}{\cm}}{\to}\quots{N}{\cm\cdot N}$, is non-surjective. Thus, for any $M\stackrel{\psi}{\to}\cm\cdot N$,
    the morphism $\phi+\psi$ remains non-surjective.
\eee

\eex

\subsection{Improving the height of an ideal}
\bprop\label{Thm.Deforming.to.expected.height}
Take a tuple of polynomials, $\uh\in R[\ux]^c$, resp.  power series, $\uh\in R[[\ux]]^c$. Assume $height(\uh)=c$.
 For any tuple $\uf\in R^r$ (for polynomials), resp.  $\uf\in \cm\cdot R^r$  (for power series), and any submodule $N\sseteq \cm\cdot R^r$ such that
 $height(Ann\quots{R^r}{N})\ge c$ holds:
\bee[1.]
\item
 There exists a tuple $\ug\in N$ such that the ideal $\big(\uh(\uf+\ug)\big)\sset R$, resp.   $\big(\uh(\uf+\ug)\big)\sset\hR^{(\cm)}$ (for power series),
   is of height $c$.
\item
 Moreover, such $\ug$ and $\uf+\ug$ can be chosen to avoid any given finite
  collection of submodules $\{N_\al\sset R^r\}$ that satisfy: $N_\al\not\supseteq \cm^n\cdot N$, for any $n,\al$.
 If $R$ is $\cm$-complete then one can avoid any countable collection of such modules.
\eee
\eprop
\bpr
The proof is similar to that of theorem \ref{Thm.Deforming.to.expected.grade}, for the case $M=R$. Recall that height does not depend on nilpotents.
 Thus at each step we pass freely to the reduced ring $\quots{R}{nilp(R)}$, and relate the grade and the height.
\bee[\rm{\bf Part}  1.]
\item
The case $c=1$. Here $h(\ux)\in R[\ux]$, $R[[\ux]]$ is a non-zero
 divisor and we should find $\ug\in N$ such that $h(\uf+\ug)\in R$ (or in $\hR^{(\cm)}$) is a non-zero divisor. As $height(Ann\quots{R^r}{N})\ge 1$, this
  annihilator contains a non-zero divisor.
   Thus (factoring by nilpotents) $grade(Ann\quots{R^r}{N})\ge 1$  and we apply Step 2.   of the proof of theorem \ref{Thm.Deforming.to.expected.grade}.

The case $c>1$. Adjust $\uf\rightsquigarrow \uf+\ug^{(1)}$ such that $h_1(\uf+\ug^{(1)})\in R$ is a non-zero divisor and pass
 to the ring $\quots{R}{h_1(\uf+\ug^{(1)})}$.
 Note that the height of the ideal $Ann\quots{R^r}{M}{\otimes} \quots{R}{h_1(\uf+\ug^{(1)})}$ is at least $c-1$, in particular this ideal contains a non-zero divisor.
  Continue as in the case $c=1$.
\item
Goes in the same way as Step 3. of theorem \ref{Thm.Deforming.to.expected.grade}. \epr
\eee

\bex\label{Ex.Deforming.to.good.height}
 We list immediate consequences, as in example \ref{Ex.Deforming.to.reg.sequence}.
 Take  $N=\oplus^r_{j=1} I_j\sset R^r$, where $I_j\sseteq \cm$ and $\{height(I_j)\ge c\}_j$.
 Then $Ann\quots{R^r}{N}=\cap I_j$ is of height $min\{height(I_j)\}$, which is at least $c$.
\bee[i.]
\item
 For  $h_1=x_1$, \dots, $h_r=x_r$, we get: for any sequence $\uf\in \cm \cdot R^r$ and any ideals $\{I_j\sset \cm \}$,
  with $\{height(I_j)\ge r\}$,
 there exist  $\{g_j\in I_j\}$ such that
 the ideal $(\uf+\ug)$ is of height $r$.

Thus being of expected height is a generic property for the ideals with a given number of generators,    in the sense of definition \ref{Def.Genericity.Property}.
\item Suppose for the sequence $f_1,\dots,f_r$, and some $k<r$ holds, $height(f_1,\dots,f_k)=k$. Then for any ideals
   $I_{k+1},\dots,I_r$ satisfying $height(I_j+\sum^k_{i=1}(f_i))\ge r$, there exist  $g_{k+1}\in I_{k+1},\dots,g_r\in I_r $ such that
\beq
height\big((f_1)+\cdots+(f_k)+(f_{k+1}+g_{k+1})+(f_r+g_r)\big)=r.
\eeq

\item
Proposition \ref{Thm.Deforming.to.expected.height} can be applied repetitively  to approximate several ideals by ideals of ``correct" height,
  as in example \ref{Ex.Deforming.to.reg.sequence.several.sequences}.
 Take some   sequences, $\{\uh^{(j)}\in R[[\ux]]^{c}\}_{j=1\dots p}$  and ideals $\{I_j\}_{j=1\dots p}$ in $R$, with $height(I_j)\ge c$.
   Fix some $r$-tuple
  $\uf\in \cm \cdot R^r$. Fix $\ug^{(1)}\in I_1\cdot R^r$ such that   $height(\uh^{(1)}(\uf+\ug^{(1)})=c$.
   Now
   shrink the ideals $I_2,\dots,I_p$ to the ideals  $I_2\cap \cm ^n,\dots,I_p\cap \cm ^n$, $n\gg1$, so that deformation by any
     $\ug^{(j)}\in I_j\cap \cm ^n$ does not decrease the height of $(\uh^{(1)}(\uf+\ug^{(1)}))$.
 (This is possible by  lemma \ref{Thm.Stability.Height.under.small.deformations}.)
     Fix $\ug^{(2)}\in (I_2\cap \cm ^n)\cdot R^r$ to ensure the needed height
     of both
 $(\uh^{(1)}(\uf+\ug^{(1)}+\ug^{(2)}))$ and  $(\uh^{(2)}(\uf+\ug^{(1)}+\ug^{(2)}))$. And so on.
\eee
\eex

\subsection{The geometric interpretation}
   Let $(R,\cm)$ be a local, complete ring, over a field $\k$.
     A tuple $\uf\in \cm\cdot R^r$
\\\parbox{12.5cm}
     {defines the morphism $\k[[\ux]]\to R$, $\ux\to \uf$. Consider the corresponding diagram of germs.
      Here the collections $\{Y_\be\}_\be$, $\{X_\al\}_\al$ are finite (or at most countable for $R$ - complete), and  these sub-germs are irreducible.
}\quad\quad $\bM Spec(R)\stackrel{\nu_\uf}{\to}Spec(\k[[\ux]])\\\cup\quad\quad\quad \cup\\\ Z,  \{Y_\be\}_\be\quad \{X_\al\}_\al\eM$

\bprop
     If $codim (Z ,o) \ge max\{codim(Y_\be)\}+max\{codim(X_\al)\}$ then
  exists
  $\ug\in   I_{(Z,o)}\cdot R^r$ such that all the intersections $\nu_{\uf+\ug}(Y_\be)\cap X_\al$ are proper,
   i.e. $codim_{Spec(R)}(Y_\be \cap \nu^{-1}_{\uf+\ug} X_\al)=codim Y_\be+codim X_\al$ for any $\al,\be$.
\eprop
 Here $codim(Z,o)$ is the minimum among the co-dimensions of the irreducible components of $(Z,o)$.
\bpr
Example     \ref{Ex.Deforming.to.good.height}.iii. reduces the question to the case of a pair of irreducible germs,
  $I_{Y}\sset R\stackrel{\uf}{\leftarrow} \k[[\ux]]\supset I_{X}$.

Fix a sequence $\uy$ in $I_Y$ such that $height(\uy)=height(I_Y)$, and a sequence $\uq$ in $I_X$ such that $height(\uq)=height(I_X)$.
 We want to find $\ug\in I_Z\cdot R^r$ such that
\beq
height(\uy,\uq(\uf+\ug))=height(I_X)+height(I_Y).
\eeq
  Define $\{h_i=y_i\}_i$ and $\{\tilde{h}_j=q_j(\ux)\}_j$, elements of $R[[\ux]]$.
    By example \ref{Ex.Deforming.to.good.height}.ii. there exist the needed $\ug$.
   Thus $codim(Y \cap \nu^{-1}_{\uf+\ug} X)=codim Y+codim X$.
\epr

\section{Applications}

\subsection{Properties of determinantal ideals}\label{Sec.Application.to.Matrices}
Let $X=\{x_{ij}\}\in Mat_{m\times n}(\k[x_{ij}])$, $m\le n$, be the matrix of indeterminates.
  Denote by $I_r(X)$ the ideal of all the $r\times r$ minors. The ideals $\{I_r(X)\}$, and the rings $\quots{\k[x_{ij}]}{I_r(X)}$
   possess numerous good properties, e.g.
\bee[i.]
\item The ideals $I_{r+1}(X)$ are prime and $grade(I_{r+1}(X))=height(I_{r+1}(X))= (m-r)(n-r)$ for $r=0,\dots,m$.
 (see e.g. Theorem 2.5 of  \cite{Bruns-Vetter}.)
\item The singular loci are: $Sing(V(I_{r+1}(X)))=V(I_r(X))\sset Spec(\k[x_{ij}])$.
\eee

 Let $(R,\cm)$ as in \S\ref{Sec.Conventions}, and  $\phi\in Mat_{m\times n}(\cm)$. In general the good properties of  $\{I_r(X)\}$
  do not survive the restriction to $\{I_r(\phi)\}$, e.g.  the grade/height equalities become inequalities:
\beq
\ber
grade_\cm(I_{j+1}(\phi))\!\le\! \min\!\!\big( (m-j)(n-j),grade_\cm(R)\big), \\
 height(I_{j+1}(\phi))\!\le\! \min\!\!\big( (m-j)(n-j),height(\cm)\big),
\\
Sing(V(I_{j+1}(\phi)))\supseteq V(I_j(\phi)).
\eer\eeq
In the degenerate case (some of) the inequalities can be strict.
 We prove that many properties are semi-continuous/stable and can be restored by deformations/approximations.

\subsubsection{Grade and height for rectangular matrices}
 Denote by $\square$ a (not necessarily square) block of some (prescribed) rows/columns, e.g. $\square$ can be a particular entry or the whole matrix.
  Thus for $\phi\in \Mat$ we have the sub-matrix $\phi_\square$, and its determinantal ideals $\{I_j(\phi_\square)\}_j$.

\bcor
For any $\phi\in Mat_{m\times n}(\cm)$ exists $q\in \N$ such that for any $\psi\in Mat_{m\times n}(\cm^q)$  holds:
\bee[1.]
\item Suppose for a set of square blocks, $\{\square_\al\}_\al$ (possibly of varying sizes), the sequence $\{det(\phi_{\square_\al})\}_\al$ is regular.
 Then the sequence $\{det(\phi_{\square_\al}+\psi_{\square_\al})\}_\al$ is regular.

\item
For any $j$ and any block $\square$ holds:
\[grade\big(I_{j+1}(\phi_\square+\psi_\square)\big)\ge grade\big(I_{j+1}(\phi_\square)\big)\quad\quad
\text{ and }
\quad\quad
height\big(I_{j+1}(\phi_\square+\psi_\square)\big)\ge height\big(I_{j+1}(\phi_\square)\big).
 \]
\eee
\ecor
Part 1 follows straight from example \ref{Ex.Preserving.Exactness.Grade.under.deformations}.i.
Part 2  follows straight from example \ref{Ex.Preserving.Exactness.Grade.under.deformations}.iv and
 lemma \ref{Thm.Stability.Height.under.small.deformations}.

\bcor\label{Thm.Det.Ideals.Improving.Grade.Height}
For any $\phi\in Mat_{m\times n}(\cm)$ and an ideal $J\sset R$  there exists  $\psi\in Mat_{m\times n}(J)$ such that $\phi+\psi$ satisfies:
\bee[1.]
\item
Let $X=\{x_{ij}\}\in Mat_{m\times n}$ be a matrix of indeterminates. Suppose for a set of square blocks $\{\square_\al\}_\al$ (possibly of varying sizes),
  the sequence
 $\{det(X_{\square_\al})\}_\al$ is regular and of size $\le grade_J(R)$. Then the sequence   $\{det(\phi_{\square_\al}+\psi_{\square_\al})\}_\al$ is regular.

\item For any $\square$ of size $\tm\times\tn$, and any $j\le min(\tm,\tn)$
holds:
\[
grade(I_{j+1}(\phi_\square+\psi_\square))= min\Big( (\tm-j)(\tn-j), grade_J(R)\Big).
\]
\item
For any $\square$ of size $\tm\times\tn$ and any $j\le min(\tm,\tn)$  holds:
\[
height(I_{j+1}(\phi_\square+\psi_\square))=min\Big( (\tm-j)(\tn-j), height(J)\Big).
\]
\eee
\ecor
Part  1 follow straight from Theorem \ref{Thm.Deforming.to.expected.grade}. (Choose $\{h_\al(\ux)=det(X_{\square_\al})\}$.)
Part  2 follow straight from example \ref{Ex.Deforming.to.reg.sequence.several.sequences}.
Part 3 follows by example \ref{Ex.Deforming.to.good.height}.

\

Sometimes one wants to deform a matrix to improve one determinantal ideal, while preserving  (some invariants of) the others. For example, we have:
\bcor
For any $\phi\in Mat_{m\times n}(\cm)$, any $r\ge1$, and any $q\in \N$, there exists $\psi\in   Mat_{m\times n}(\cm^q\cdot I_r(\phi))$ such that:
\bee[i.]
\item $ I_j(\phi+\psi)=  I_j(\phi)$ for any  $j=1,\dots,r$.
\item $ grade\big( I_{j+1}(\phi+\psi)\big)= \min\big((m-j)(n-j),grade(I_r(\phi))\big)$  for any $j\ge r$.
\item  $height\big( I_{j+1}(\phi+\psi)\big)= \min\big((m-j)(n-j),height(I_r(\phi))\big)$ for any $j\ge r$.
\eee
\ecor
(Choose $J=I_r(\phi)$ in the last Corollary.)

\subsubsection{Regularity properties}
\bcor
Let $(R,\cm)$ be a local Cohen-Macaulay ring,  $dim(R)\ge (m-r)(n-r)$.
 Then for any $\phi\in Mat_{m\times n}(\cm)$  and any $q\in \N$ there exists $\psi\in Mat_{m\times n}(\cm^q\cdot I_r(\phi))$ such that all
  the ideals $\{I_{j+1}(\phi+\psi)\}_{j\ge r}$ and the rings $\{\quots{R}{I_{j+1}(\phi+\psi)}\}_{j\ge r}$ are perfect Cohen-Macaulay.
\ecor
\bpr
Apply the Eagon-Hochster criterion: if $R$ is Cohen-Macaulay and $height(I_{r+1}(\phi))=(m-r)(n-r)$ then $I_{r+1}(\phi)$ and $\quots{R}{I_{r+1}(\phi)}$
 are Cohen-Macaulay.

Similarly, $I_r(\phi)$ is perfect if it has the expected grade.
\epr

\subsubsection{The case of (skew-)symmetric matrices}
 In this case the statements are similar to the corollaries above, we just recall the heights for matrices of indeterminates:
\bee[i.]
\item
For any $X=\{x_{ij}\}\in Mat^{sym}_{m\times m}(\k[x_{ij}])$
 the ideals $I_{j+1}(X)$ are prime and $height(I_{j+1}(X))=\bin{m+1-j}{2}$ for $j=0,\dots,m$.

\item
For any $X=\{x_{ij}\}\in Mat^{skew-sym}_{m\times m}(\k[x_{ij}])$:
 \bei
 \item ($j$ - even) $height(I_{j+1}(X))=\bin{m-j}{2}$.
 \item ($j$ - odd) $height(I_{j+1}(X))=\bin{m+1-j}{2}$.
 \eei
\eee

\subsection{The grade of a module and height of an ideal}\label{Sec.Depth.Grade.Improving}
Let $(R,\cm)$ be a local ring,
    take an ideal  $I\sset R $, a module $M$ and a sequence $\uf\in \cm \cdot R^r$.
\bee[\bf i.]
\item
By \S\ref{Sec.Depth.Grade.Semi-continuity} there exists $n$ such that for any $\ug\in \cm^n\cdot R^r$ holds: $grade_{(\uf+\ug)}M\ge grade_{(\uf)}M$.

Obviously $grade_{(\uf)}M\le \min( r,grade_\cm M)$ and the strict inequality occurs in degenerate cases.
  However, for any $J\sset \cm\cdot R^r$ such that $grade_J M=\min( r,grade_\cm M)$,
  there exists a tuple $\ug\in J\cdot R^r$ such that $grade_{(\uf+\ug)}M=min(r,grade_\cm M)$.
 (Apply example \ref{Ex.Deforming.to.reg.sequence}.)

\item
Take a sequence $\uf\in \cm\cdot R^r$ and the corresponding projection $R\stackrel{\pi_\uf}{\to}\quots{R}{(\uf)}$.
 Assume $grade_I M\ge r$ then $grade_{\pi_\uf(I)}\pi_\uf M\ge grade_I M-r$.
  For generic $\uf$ one expects the equality.
     We claim:
 if $ grade_I R+r\le grade_\cm R$ then
     for any submodule $N\sset \cm\cdot R^r$, such that $Ann\quots{R^r}{N}$ contains $M$-regular sequence of length $grade_I M$,
     exist $\ug\in N$ such that $grade_{\pi_{\uf+\ug}(I)}\pi_{\uf+\ug} M= grade_I M-r$.

\bpr Let $grade_I M=c\ge r$. Take an $M$-regular sequence $q_1,\dots,q_{c-r}\in I$. The sequence $\uq,\uf$
 is not necessarily $M$-regular. By example \ref{Ex.Deforming.to.reg.sequence} there exists $\ug\in N$ such that
   $\uq,\uf+\ug$ is $M$-regular and also the sequence $\uf+\ug$ is $\quots{R}{I}$-regular. Thus
 $grade_{\pi_{\uf+\ug}(I)}\pi_{\uf+\ug} M\ge grade_I M-r$.

 The other direction follows from the relation: $grade_I \quots{M}{f\cdot M}=grade_I M-1$, for $f$ generic.\epr


\item
 Suppose $height(I)\ge r$.
   Then  $   height(\pi_\uf I)\ge height(I)-r$.
 For generic $\uf$ one expects the equality. Proposition \ref{Thm.Deforming.to.expected.height} gives:
  for any submodule $N\sset \cm\cdot R^r$, such that $height(Ann\quots{R^r}{N})\ge r$,       exists $\ug\in N$ such that
  $   height(\pi_{\uf+\ug}I)= height(I)-r$.

Proof: as above.
\eee

\subsection{Properties of morphisms  and generalized Eagon-Northcott complexes}\label{Sec.Application.to.E.N.complexes}
To a morphism of free modules,   $F_1\stackrel{\phi}{\to}F_0$, $rank(F_1)\ge rank(F_0)$,
 one associates the generalized Eagon-Northcott complexes, $(\cC^i_\phi)_\bullet$, for $i=0,1,\dots,rank(F_0)$, see Appendix $A.2.6$ of  \cite{Eisenbud}.
 In particular, $(\cC^0_\phi)_\bullet$ is the E.N. complex (in the generic case it resolves the ideal $I_m(\phi)$), while $(\cC^1_\phi)_\bullet$
  is the Buchsbaum-Rim complex (in the generic case it resolves $coker(\phi)$).

In non-generic cases the complexes $\cC^i_\phi$  are not necessarily exact. However, this can be improved.
\bcor
\bee[1.]
\item
 (Stability of exactness) Suppose the complex $\cC^i_\phi$ is exact, for some $i$. Then exists $n$ such that for any  $F_1\stackrel{\psi}{\to}\cm^n\cdot F_0$
 the complex $\cC^i_{\phi+\psi}$ is exact.
 \item (Exact approximation) Suppose     $grade_\cm(R)\ge  (rank(F_1)-r)(rank(F_0)-r)$, for some $r$.
 Then for any $n\in \N$  there exists $F_1\stackrel{\psi}{\to}\cm^n\cdot F_0$ such that the  complexes $\{\cC^i_{\phi+\psi}\}_{i> r}$ are acyclic.
\eee
\ecor
\bpr
{\bf 1.}
The maps in   $\cC^i_\phi$  are constructed via the exterior powers of $\phi$. Thus a deformation $\phi\rightsquigarrow \phi+\psi$,
   induces the deformed complexes $\cC^i_{\phi+\psi}$. Apply lemma \ref{Thm.Stability.of.exactness.under.small.deformations}.

{\bf 2.}
By Corollary \ref{Thm.Det.Ideals.Improving.Grade.Height}, exists $F_1\stackrel{\psi}{\to}J\cdot F_0$ such that for any $i\ge r$ holds:
\beq
grade_{I_{i+1}(\phi+\psi)}R=(rank(F_1)-i)(rank(F_0)-i).
\eeq
 And this ensures the acyclicity of $\{\cC^i_{\phi+\psi}\}_{i> r}$, see Theorem A.2.10 of \cite{Eisenbud}.\epr

\bex
Let $M$ be a module with the Betti numbers $\be_1(M)\ge \be_0(M)$, and assume $grade_\cm R\ge \be_1-\be_0+1$. Then $M$
is approximated by a module $\tM$ with  $\be_i(\tM)=\be_i(M)$, for $i=0,1$, and the projective dimension  $pd(\tM)=\be_1-\be_0+1$.
\eex

\subsection{ $Tor,Ext$}\label{Sec.Application.to.Tor.Ext}

\bex
 (Stability of $Tor$, $Ext$-vanishing) Take a morphism of free modules,   $F_1\stackrel{\phi}{\to}F_0$, and assume  (for some $i$) the generalized  E.N. complex
 $\cC^i_\phi$ is acyclic.
 \bee[\bf i.]
 \item  Assume $Tor_j^R(coker(\phi),M)=0$ for some $j\in\N$.
 Then there exists $n\in \N$ such that for any  $F_1\stackrel{\psi}{\to}\cm^n\cdot F_0$ holds: $Tor_j^R(coker(\phi+\psi),M)=0$ (for those $j$ and $M$).

Indeed, resolve $coker(\phi)$ by  $ \cC^i_\phi $, then the complex $ \cC^i_\phi \otimes M$ is exact at place $j$. Thus any $\cm^n$-adic approximation,
 $\cC^i_{\phi+\psi}\otimes M$ is exact at place $j$.
\item
Similarly one gets: if the complex  $ \cC^i_\phi $ is acyclic and $Ext_R^j(coker(\phi),M)=0$, for some $j,M$, then $Ext_R^j(coker(\phi+\psi),M)=0$
 for any  $F_1\stackrel{\psi}{\to}\cm^n\cdot F_0$.
\item
In the particular case, $rank(F_0)=1$, $rank(F_1)=r$, we get again: if $\uf$ is $R$-regular and $M$-regular, then $\uf+\ug$ is $R$ and $M$-regular.
\item
The acyclicity assumption on $ \cC^i_\phi $ is   important due to the following trivial example. Let $R\stackrel{\phi=0}{\to}R$ and
 suppose $\psi\in \cm^n$ is regular. Then $Tor_1(coker(\phi),M)=0$, but $Tor_1(coker(\phi+\psi),M)=0:_M\psi$.
\eee\eex

\bex (Improving $Tor$ and $Ext$)
\bee[\bf i.]
\item
 Theorem \ref{Thm.Deforming.to.expected.grade} implies: if $r\le grade_\cm(R), grade_\cm(M)$
  then for any $\uf\in \cm\cdot R^r$ and $n\in \N$ exists $\ug\in \cm^n\cdot R^r$ such that the sequence $\uf+\ug$ is $R$-regular and $M$-regular.

Therefore the functor $ \quots{R}{(\uf+\ug)}\otimes\underline{\ \ }$ preserves exactness of the minimal resolution of $M$.
 Hence $Tor^R_i(\quots{R}{(\uf+\ug)},M)=0$ for any $i>0$.
\item Suppose $grade_\cm M\ge r$. Then for any $\uf\in \cm\cdot R^r$ and $n\in \N$ exists $\ug\in \cm^n\cdot R^r$
 such that
 \beq
 Ext^i_R(\quots{R}{(\uf+\ug)},M)=0 \quad \text{ for }\quad i<r.
 \eeq
  Indeed, by \S\ref{Sec.Depth.Grade.Improving}
  there exists $\ug\in \cm^n\cdot R^r$ such that $grade_{(\uf+\ug)}M=r$. It suffices to recall an equivalent definition,
  $grade_I M:=inf\{i|\  Ext^i_R(\quots{R}{I},M)\neq0\}$.

\eee
\eex

\end{document}